\documentclass[11pt,a4paper]{article}
\usepackage[a4paper,margin=1in]{geometry}
\usepackage{ifthen,latexsym,amssymb,amsmath,bbm,fixmath,stmaryrd}
\usepackage[backend=bibtex, style=numeric, giveninits=true]{biblatex}
\usepackage[shortlabels]{enumitem} 
\usepackage{listings}
\usepackage{xcolor}

\definecolor{codegreen}{rgb}{0,0.6,0}
\definecolor{codegray}{rgb}{0.5,0.5,0.5}
\definecolor{codepurple}{rgb}{0.58,0,0.82}
\definecolor{backcolour}{rgb}{0.95,0.95,0.92}

\lstdefinestyle{mystyle}{
    backgroundcolor=\color{backcolour},   
    commentstyle=\color{codegreen},
    keywordstyle=\color{magenta},
    numberstyle=\tiny\color{codegray},
    stringstyle=\color{codepurple},
    basicstyle=\ttfamily\footnotesize,
    breakatwhitespace=false,
    breaklines=true,
    keepspaces=true,
    showspaces=false,
    showstringspaces=false,
    showtabs=false,
    tabsize=4
}

\usepackage{algorithm2e}
\usepackage{tabularx}
\usepackage{float}
\usepackage{etoolbox}
\usepackage{changepage}

\usepackage{pdflscape}
\usepackage{hyperref}
\hypersetup{
  colorlinks=true,
  linkcolor=blue,
  citecolor=blue,
  urlcolor=blue
}
\usepackage[noabbrev]{cleveref}

\begin{document}

\newcommand{\hide}[1]{}

\title{\texttt{FlagAlgebraToolbox}: Flag Algebra Computations in SageMath}

\author{Levente Bodnár}

\maketitle

\begin{abstract}
We introduce \texttt{FlagAlgebraToolbox}, an extension of SageMath capable of automating flag algebra calculations and optimizations. \texttt{FlagAlgebraToolbox} has a simple interface, can handle a wide range of combinatorial theories, can numerically optimize extremal combinatorial problems and round the results to produce exact proofs. We present the core concepts used in the toolbox, with example workflows. 
\end{abstract}

\section{Introduction}\label{sec:introduction}

Flag algebras, introduced by Razborov~\cite{Razborov07}, provide a formal calculus describing the relations between densities of small patterns in limit combinatorial structures. This framework can turn extremal combinatorial questions into semidefinite optimization problems providing rigorous inequalities and often sharp asymptotic bounds. In practice, carrying out a flag algebra computation involves a large number of small steps suitable for computer automation: 
\begin{itemize}
    \item generating a list of combinatorial objects and flags,
    \item calculating the chain rule, multiplication relations and projections between the flags,
    \item finding optimal sum-of-squares expressions.
\end{itemize}
These steps are conceptually uniform across many settings, yet existing implementations \cite{flagmatic} are often tailored to a small number of standard theories (such as graphs or $k$-graphs).

This paper describes \texttt{FlagAlgebraToolbox}, an extension of SageMath~\cite{The_SageMath_Developers_SageMath_2026}, designed to support the flag algebra formal calculus in an extensible way. The toolbox provides high-level methods for constructing combinatorial theories (including graphs, digraphs, hypergraphs, hypercubes), generating and manipulating flags and typed expressions symbolically. The package is designed to support major techniques used in optimizations involving flag algebras. 

\subsection{Outline}\label{subsec:outline}
The next subsection~\ref{subsec:installation} describes the recommended installation on Linux, macOS, and Windows, as well as the source distribution. Subsection~\ref{subsec:quick_start} showcases the standard workflow, setting up and solving Mantel's problem~\cite{Mantel1907}, and an example problem on 3-graphs. Section~\ref{sec:concepts} introduces the core abstractions used throughout the package: flags and types, combinatorial theories, and tools to solve combinatorial optimization problems, including rounding, constructions and certificates. Section~\ref{sec:examples} presents illustrative examples and provides scripts used in recent work. We conclude with limitations and directions for future work in Section~\ref{sec:future}.

% Replace the complete Installation subsection with the following.
\subsection{Installation}\label{subsec:installation}

The recommended installation is the prebuilt Docker image, providing a SageMath environment on the major desktop operating systems including the flag algebra toolbox functionality. The launcher script
\texttt{start-fat.sh} is available at \url{https://gist.github.com/bodnalev/900b7f8c1b4dcf0e69e232be361c04ff}.

\paragraph{Linux and macOS}
Install Docker Engine or Docker Desktop and ensure that Docker is running; see \url{https://docs.docker.com/get-started/get-docker/}. Download the launcher under the filename \texttt{start-fat.sh}, make it executable, and run it from a terminal:
\begin{lstlisting}[language=bash]
chmod +x start-fat.sh
./start-fat.sh
\end{lstlisting}
On some Linux versions, access to the Docker thread initially requires sudo privileges. Users can allow their account to use Docker without \texttt{sudo} using:
\begin{lstlisting}[language=bash]
sudo usermod -aG docker "$USER"
newgrp docker
\end{lstlisting}

\paragraph{Windows}
On Windows, the recommended setup is Docker Desktop using its WSL~2 backend, together with a WSL~2 Linux distribution in which the Bash launcher is run. To install WSL, from an administrator PowerShell window run
\begin{lstlisting}
wsl --install
\end{lstlisting}
Start Docker Desktop and enable integration with the installed WSL distribution under \texttt{Settings > Resources > WSL Integration}. The launcher can then be run from the Ubuntu/WSL terminal exactly as on Linux. The URL printed by the launcher can be opened in a normal Windows browser. For best file-system performance, it is preferable to keep the notebook directory inside the WSL Linux file system rather than under \texttt{/mnt/c}. Further instructions are available at \url{https://docs.docker.com/desktop/features/wsl/} and \url{https://learn.microsoft.com/windows/wsl/install}.

\paragraph{Using the launcher}
On its first run, the launcher pulls the docker image. Passing \texttt{--update} attempts to download an update if available:
\begin{lstlisting}[language=bash]
./start-fat.sh --update
\end{lstlisting}
The launcher starts Jupyter, and normally opens it in the default browser. To start without the browser, use the option \texttt{--no-browser}. The terminal should be kept open while Jupyter is in use; pressing \texttt{Ctrl+C} stops Jupyter and the container.

\paragraph{Building from source}
The source is available at \url{https://github.com/bodnalev/sage}. Building the SageMath fork requires a full compiler toolchain, substantial disk space, and considerably more time than using the container image. Current source-build requirements and commands are available in the repository and in the SageMath installation guide at \url{https://doc.sagemath.org/html/en/installation/source.html}.

\subsection{Quick start}\label{subsec:quick_start}

Running the launcher opens a local Jupyter file explorer rooted at the selected notebook directory. Creating a new notebook with the SageMath kernel starts an interactive environment. Running following script calculates that the maximum density of edges in triangle-free graphs is at most $1/2$.
\begin{lstlisting}[language=Python]
triangle = GraphTheory(3, edges=[[0,1],[0,2],[1,2]])
edge = GraphTheory(2, edges=[[0,1]])
GraphTheory.exclude(triangle)
bound = GraphTheory.optimize(edge, 3, maximize=True, exact=True)
print(bound)
GraphTheory.reset()
\end{lstlisting}
The first two lines define the triangle and edge graphs. \texttt{GraphTheory} is the theory object for simple undirected graphs. This theory contains one binary relation called \texttt{edges} which is symmetric and irreflexive. To define a \texttt{GraphTheory} object, we can call the theory and provide the vertex number and the pairs forming the \texttt{edges} relation. If the vertex number $N$ is provided, then the relations expect elements from the vertex set $\{0, 1, ..., N-1\}$. In this example, the variables \texttt{triangle} and \texttt{edge} store the corresponding triangle and edge graphs.

The third line modifies the theory to only contain graphs without triangles, by excluding the \texttt{triangle} object. The final line, in this theory without triangles, tries to maximize the density of edges. The second parameter $3$ in the function indicates that the optimization is performed by truncating the flag algebra to graphs with vertex number at most $3$. By default the optimization is only numeric. Here, we apply an automatic rounding method with the parameter \texttt{exact=True}, to get an exact result. Finally the found bound $1/2$ printed. The last line resets the theory, so the triangle is no longer excluded (see subsection~\ref{subsec:theories} for more details).

In \texttt{FlagAlgebraToolbox}, the combinatorial objects are induced. The following example maximizes in 3-graphs the induced density of 4-tuples containing exactly one edge, under the condition that the density of edges is at least $1/2$.
\begin{lstlisting}[language=Python]
target = ThreeGraphTheory(4, edges=[[0,1,2]])
edge = ThreeGraphTheory(3, edges=[[0,1,2]])
expression = edge - 1/2
ThreeGraphTheory.optimize(target, 6, maximize=True, positives=[expression])
\end{lstlisting}
Similarly to \texttt{GraphTheory}, the \texttt{ThreeGraphTheory} encodes a theory for $3$-graphs. There is a single ternary relation which is symmetric and irreflexive in a sense that double and triple entries are not allowed. The first line defines the target $3$-graph structure, with $4$ vertices and one edge. The second line defines the $3$-edge. The third line defines an expression which is the difference of an \texttt{edge} and the constant $1/2$. The final line runs the optimization, with the \texttt{positives=[expression]} parameter indicating the constraint that the edge density must be at least $1/2$, i.e. the \texttt{expression} must evaluate to a non-negative expression (see subsection~\ref{subsec:optimization} for more details).

Under these constraints, the optimizer returns the numeric value approximately \texttt{0.563219049}, which is a numeric upper bound. With the parameter \texttt{exact=True}, it gives a weak rational upper bound, which is $8701/15360$, approximately \texttt{0.566471354}.

\section{Core concepts}\label{sec:concepts}

This section describes the main building blocks and abstractions used in the package.  We separate the combinatorial layer (theories, flags, types, patterns) from the algebraic layer (flag algebra elements over a chosen base field) and the optimization layer (formulating and solving SDP-based extremal problems).  Before the detailed description of the theories in subsection~\ref{subsec:theories}, the examples will use \texttt{GraphTheory}. In particular, they will use the shortcut
\begin{lstlisting}[language=Python]
G = GraphTheory
\end{lstlisting}
to construct these objects with a shorter syntax.

\subsection{Working with flags}\label{subsec:flags}

\paragraph{Creating flags}
A flag is created by specifying the size together with the relevant relations. To create a flag in a given \texttt{Theory}, the theory object is called \texttt{Theory(size, **relations)} with the desired size and the relations as named parameters. 
\begin{lstlisting}[language=Python]
# graphs
triangle = G(3, edges=[[0,1],[0,2],[1,2]])
cherry = G(3, edges=[[0,1],[0,2]])
# 3-graphs
k4m = ThreeGraphTheory(4, edges=[[0,1,2],[0,1,3],[1,2,3]])
\end{lstlisting}

Flags are considered up to isomorphism inside the underlying theory; in particular, relabeling vertices does not change the represented element. Vertices of a flag can be marked, forming the flag's type. A type is defined with the \texttt{ftype} parameter. The order of the marked vertices matters. The equality respects the types. In the example below, the two different pointed versions of the cherry graph are different.  
\begin{lstlisting}[language=Python]
pointed_edge = G(2, edges=[[0,1]], ftype=[0])
pointed_cherry = G(3, edges=[[0,1],[1,2]], ftype=[1])
other_pointed_cherry = G(3, edges=[[0,1],[1,2]], ftype=[0])
\end{lstlisting}
A flag in which all vertices are marked is simply a type. We can get the type of any flag by calling the \texttt{.ftype()} function. 
\begin{lstlisting}[language=Python]
G(1, ftype=[0])
pointed_cherry.ftype()
other_pointed_cherry.ftype()
\end{lstlisting}
all return the same type formed from a single point.

All flags are induced by default. Patterns allow optional relations.  A pattern can be expanded into the list or sum of all induced flags compatible with it. Patterns are defined within a theory with the \texttt{Theory.pattern(size, **relations)} syntax, where the relations are assumed to be optional by default. 
\begin{lstlisting}[language=Python]
# induced, all other edges are absent
cherry = G(3, edges=[[0,1],[1,2]]) 
# non-induced: other edges are unspecified
cherry_pattern  = G.pattern(3, edges=[[0,1],[1,2]]) 

# list compatible induced flags (cherry, triangle)
cherry_pattern.compatible_flags() 
\end{lstlisting}

Patterns may also enforce missing relations explicitly, using the suffix
\texttt{\_missing} (or \texttt{\_m} for short). The example below defines all $3$ vertex graphs at least one edge and at least one missing edge. So it matches the cherry and the cherry complement.
\begin{lstlisting}[language=Python]
tpat = G.pattern(3, edges=[[0,1]], edges_missing=[[1,2]])
tpat.compatible_flags()
\end{lstlisting}

\subsection{Flag algebras and flag algebra elements}
Flags represent combinatorial structures. However, they often get automatically coerced into a \texttt{FlagAlgebraElement} to perform arithmetic on them. 

As an example, consider the sum of an edge and a triangle. This results in a linear combination of flags, which is a \texttt{FlagAlgebraElement} living in a suitable \texttt{FlagAlgebra} parent object. 

\begin{lstlisting}[language=Python]
edge = G(2, edges=[[0,1]]) 
triangle = G(3, edges=[[0,1],[1,2],[2,0]])
print(edge+triangle)
\end{lstlisting}
\begin{lstlisting}
Flag Algebra Element over Rational Field
0   - Flag on 3 points, ftype from () with edges=()
1/3 - Flag on 3 points, ftype from () with edges=(01)
2/3 - Flag on 3 points, ftype from () with edges=(01 02)
2   - Flag on 3 points, ftype from () with edges=(01 02 12)
\end{lstlisting}

The resulting flag algebra element is expressed as a sum of $3$ vertex flags. Standard arithmetic operations work on flags with the same type. The following example calculates the square of a pointed edge.
\begin{lstlisting}[language=Python]
pointed_edge = G(2, edges=[[0,1]], ftype=[0])
print(pointed_edge * pointed_edge)
\end{lstlisting}
\begin{lstlisting}
Flag Algebra Element over Rational Field
0 - Flag on 3 points, ftype from (0,) with edges=()
0 - Flag on 3 points, ftype from (0,) with edges=(01)
0 - Flag on 3 points, ftype from (2,) with edges=(01)
1 - Flag on 3 points, ftype from (0,) with edges=(01 02)
0 - Flag on 3 points, ftype from (1,) with edges=(01 02)
1 - Flag on 3 points, ftype from (0,) with edges=(01 02 12)
\end{lstlisting}

The projection of a flag with a large type to a smaller one is done by the \texttt{project()} function. Calling it without parameters projects to the empty type. The following example proves Mantel's theorem with these basic operations.
\begin{lstlisting}[language=Python]
triangle = G(3, edges=[[0,1],[1,2],[2,0]])
edge = G(2, edges=[[0,1]])
G.exclude(triangle)
pointed_edge = G(2, edges=[[0,1]], ftype=[0])
degree_imbalance = pointed_edge - 1/2
print(edge - 1/2 + 2 * (degree_imbalance * degree_imbalance).project())
G.reset()
\end{lstlisting}
This script calculates the square of the degree imbalance flag algebra element, which simply counts at a marked point how much the relative degree differs from $1/2$. The projection of this square is guaranteed to be positive, and provides a sufficient lower bound to the claim. The printed flag algebra element has all coefficient negative arguing that the density of edges can not be larger than $1/2$.

The suitable algebra is automatically found to be the flag algebra constructed from linear combinations of the flags with the appropriate type. These algebras can be explicitly constructed in the package as \texttt{FlagAlgebra(BaseTheory, BaseRing, ftype)}, but these are usually found automatically during the calculation. The base ring must contain the rationals.

The following example demonstrates a calculation with the polynomial ring \texttt{R.<x> = QQ[]} as base.
\begin{lstlisting}[language=Python]
R.<x> = QQ[]
pointed_edge = G(2, edges=[[0,1]], ftype=[0])
triangle = G(3, edges=[[0,1],[1,2],[2,0]])
print(((pointed_edge-x)*(pointed_edge-1/2+x)).project() + triangle)
\end{lstlisting}
\begin{lstlisting}
Flag Algebra Element over Univariate Polynomial Ring in x over Rational Field
-x^2 + 1/2*x       - Flag on 3 points, ftype from () with edges=()
-x^2 + 1/2*x - 1/6 - Flag on 3 points, ftype from () with edges=(01)
-x^2 + 1/2*x       - Flag on 3 points, ftype from () with edges=(01 02)
-x^2 + 1/2*x + 3/2 - Flag on 3 points, ftype from () with edges=(01 02 12)
\end{lstlisting}

\subsection{Theories}\label{subsec:theories}

A theory specifies the combinatorial setting, where the combinatorial calculations are performed. This includes relation names, arities, and symmetry/ordering conventions, together with any global restrictions such as forbidden configurations. All flags, patterns, generation, and optimization routines are parameterized by the current theory state.

\paragraph{Creating theories}
Several common theories are provided out of the box, including (di)graphs and $k$-graphs, as well as vertex-color theories (unary relations) intended for colored and multipartite settings:
\begin{lstlisting}[language=Python]
GraphTheory = Theory(
    "Graph", relation_name="edges", 
    arity=2, is_ordered=False
)
DiGraphTheory = Theory(
    "DiGraph", relation_name="edges", 
    arity=2, is_ordered=True
)
ThreeGraphTheory = Theory(
    "ThreeGraph", relation_name="edges", 
    arity=3, is_ordered=False
)
Color0 = Theory("Color0", relation_name="C0", arity=1)
Color1 = Theory("Color1", relation_name="C1", arity=1)
\end{lstlisting}

New theories can be created by specifying a name, a relation name, arity, and whether the relation is ordered. For example:
\begin{lstlisting}[language=Python]
T0 = Theory("OtherDiGraph", arity=2, is_ordered=True,relation_name="diedge")
T1 = Theory("OtherThreeGraph", arity=3, relation_name="edges3")
\end{lstlisting}

There is a simple method to combine multiple theories and their relations together. This requires theories with different relation names. For example, one can combine a graph relation \texttt{edges} with a 3-uniform hypergraph relation \texttt{edges3} from the previous example.
\begin{lstlisting}[language=Python]
T0 = Theory("OtherThreeGraph", arity=3, relation_name="edges3")
CombinedTheory = combine("TwoThreeGraph", GraphTheory, T0)
\end{lstlisting}

For combined theories, flags (and patterns) can specify each relation separately:
\begin{lstlisting}[language=Python]
test_flag = CombinedTheory(3, edges=[[0,1],[0,2]], edges3=[[0,1,2]])
test_pat  = CombinedTheory.pattern(4, edges=[[0,1]], edges3_m=[[1,2,3]])
\end{lstlisting}

When combining relations of the same arity, one may optionally impose symmetries (e.g.\ full symmetry or cyclic symmetry) to identify equivalent relation-labellings and reduce the count of non-isomorphic structures:
\begin{lstlisting}[language=Python]
G = GraphTheory
Gp = Theory(
    "OtherGraph", relation_name="oedges", 
    arity=2, is_ordered=False
)
C0 = combine("DoubleEdgeGraph", G, Gp, symmetries=NoSymmetry)
C1 = combine("SymmetricDoubleEdgeGraph", G, Gp, symmetries=FullSymmetry)
\end{lstlisting}

In this case, \texttt{C0(2, edges=[[0,1]], oedges=[])} and \texttt{C0(2, edges=[], oedges=[[0,1]])} for example, are distinguished, but \texttt{C1(2, edges=[[0,1]], oedges=[])} and \texttt{C1(2, edges=[], oedges=[[0,1]])} are equal.

\paragraph{Excluding and generating flags}

Theories support forbidding induced configurations with an empty type. This is done in-place and persists until reset. Exclusions modify the theory globally; use \texttt{Theory.reset()} to restore to defaults. This forces all calculations to happen in this restricted theory, raising an error when needed. Exclusion of multiple elements is supported, by passing an iterator, this can contain patterns. For patterns, any compatible flag is excluded.  The method \texttt{generate(n, ftype)} lists all (induced) flags of order $n$ in the theory with the matching type, respecting any current exclusions. If the type is not provided, it is assumed to be the empty type.

\begin{lstlisting}[language=Python]
test_pattern = G.pattern(4, edges=[[0,1],[0,2]], edges_m=[[0,3]])
G.exclude(test_pattern)
flag_list = G.generate(5)
print("\n".join(map(str, flag_list)))
G.reset()
\end{lstlisting}

\begin{lstlisting}
Flag on 5 points, ftype from () with edges=()
Flag on 5 points, ftype from () with edges=(01)
Flag on 5 points, ftype from () with edges=(02 14)
Flag on 5 points, ftype from () with edges=(01 02 03 04)
Flag on 5 points, ftype from () with edges=(01 02 03 04 12 13 14 23 24 34)
\end{lstlisting}

\subsection{Optimization}\label{subsec:optimization}

At a high level, extremal problems are specified by a target expression and a target size, and then solved by forming the corresponding SDP.  The main function for this is \texttt{Theory.optimize}. The target may be any linear combination of flags:
\begin{lstlisting}[language=Python]
target = G(2, edges=[[0,1]]) + G(3, edges=[[0,1]])
G.optimize(target, 4)
\end{lstlisting}
This returns a numerical bound around $1.25000000116$, which is accurate up to the accumulated floating point errors. The exact optimum $5/4$ is attained at the union of two equal cliques. The \texttt{Theory.optimize} function performs a maximization by default, minimization is done by the parameter \texttt{maximize=False}.

When \texttt{exact=True} is provided, the optimizer attempts to round a numerical solution to an exact certificate, a rational number by default. For simple problems, this can find the exact optimal solution automatically.
\begin{lstlisting}[language=Python]
target = G(2, edges=[[0,1]]) + G(3, edges=[[0,1]])
G.optimize(target, 4, exact=True)
\end{lstlisting}
The above code correctly finds the exact optimal upper bound $5/4$. This is done by attempting to convert the numerical solution into a verified certificate by rounding all numeric data (scalars and semidefinite matrices) to exact rational values, while preserving the constraints in exact arithmetic. The returned bound is always verified. In simple problems the rounded result is best possible, however in complex cases this bound might be worse than the numeric value. The parameter \texttt{denom} controls the maximum denominator used in the rational approximation during rounding, increasing this value can help improve the gap between the numeric and the exact bound. 
\begin{lstlisting}[language=Python]
p4 = G(4, edges=[[0,1],[1,2],[2,3]])
bd_num = G.optimize(p4, 6)
print("Numeric bound is {}".format(bd_num))
bd_weak = G.optimize(p4, 6, exact=True, denom=512)
print("Weaker bound is {}~{}".format(bd_weak, bd_weak.n()))
bd_strong = G.optimize(p4, 6, exact=True, denom=1048576)
print("Stronger bound is {}~{}".format(bd_strong, bd_strong.n()))
\end{lstlisting}
\begin{lstlisting}
Numeric bound is 0.21357246423780124
Weaker bound is 111/512~0.216796875000000
Stronger bound is 223949/1048576~0.213574409484863
\end{lstlisting}

The inducibility of the $4$-vertex path is an interesting open problem. The best construction comes from Even-Zohar and Linial \cite{Even_Zohar_2014} giving around $0.2014$. Considering the flag algebra truncated to $8$ vertex graphs (instead of the example script above with $6$) gives the tighter upper bound $0.2045$ matching the one obtained by Flagmatic \cite{flagmatic}.

\paragraph{Assumptions} Additional assumptions can be supplied as a list of flag algebra elements, during the optimization these are assumed to be positive. For example, passing the constraint \texttt{[1/2 - edge]} to the parameter \texttt{positives} enforces edge density at most $1/2$:
\begin{lstlisting}[language=Python]
G.optimize(triangle, 4, positives=[1/2 - edge])
\end{lstlisting}
This code returns numerically $0.3535533923$, with the true optimum being $1/\sqrt{8}$ at a clique with relative size $1/\sqrt{2}$ following from a corollary of the Kruskal-Katona theorem.

The positivity constraints expect a list, adding multiple constraints optimizes under the assumption that they are all satisfied. Adding a typed expression in the constraint list translates roughly to the condition that all embeddings of the type in the large structure satisfies the property. As an example, consider the following:
\begin{lstlisting}[language=Python]
pointed_edge = G(2, edges=[[0,1]], ftype=[0])
G.optimize(triangle, 6, positives=[1/2 - pointed_edge, 1/3 - edge])
\end{lstlisting}
This maximizes the triangles, under the constraint that:
\begin{itemize}
    \item \texttt{1/2-pointed\_edge}, which translates to every vertex (the point type) having relative degree at most $1/2$,
    \item \texttt{1/3-edge}, which translates to having at most $1/3$ edge density.
\end{itemize}
The optimizer returns around $0.1527374406$, but it is an interesting question what the optimum might be.

\paragraph{Constructions}
The rounding can be further helped by providing a construction achieving the optimum. Constructions can be defined using the \texttt{Theory.blowup\_construction} function, which returns blow-ups of small templates. The syntax is similar to the flag generation. As an example, the script below constructs the $4$ vertex description of a graphon with $3$ equal parts (simply indexed by $\{0, 1, 2\}$) and the edge relations between them:
\begin{lstlisting}[language=Python]
constr = G.blowup_construction(4, 3, edges=[[0,0],[1,2]])
\end{lstlisting}
Notice that the relations here may not be irreflexive, \texttt{[0,0]} indicates that part $0$ contains a clique. The first parameter represents the size of the returned flag algebra element. In particular, \texttt{constr} defined above is a linear combination of $4$ vertex graph flags. The second parameter defines the parts of the limit structure. When a number is provided, it assumes equal sized parts, when a list is provided the part sizes are defined by the numbers in the list. 

When a construction is defined, it can be used in an optimization, to help the rounding. This is only beneficial when the numeric bound is the same as to the value provided by the construction up to machine error, usually a $10^{-8}$ gap is required. As an example, consider the problem of maximizing triangles with $1/2$ bounded edge density. In this case, providing the construction helps the rounding:

\begin{lstlisting}[language=Python]
sq2 = QQ[sqrt(2)](sqrt(2))
constr = G.blowup_construction(4, [1/sq2, 1-1/sq2], edges=[[0,0]])
G.optimize(triangle, 4, positives=[1/2-edge], 
           construction=constr, exact=True
          )
\end{lstlisting}
This successfully returns \texttt{1/4*sqrt2}, as expected. Note that the size of the construction (not the number of parts) must match the size used in the optimization. In the example above, size $4$ was used.

Quasi-random constructions can be defined by replacing the relation list, with a relation dictionary, where the keys are relation pairs, and the values are the probabilities. The following example creates a quasi-random bipartite graph, where the probability of an edge between two vertices in opposite parts is $1/3$.
\begin{lstlisting}[language=Python]
G.blowup_construction(4, 2, edges={(0, 1): 1/3})
\end{lstlisting}

Constructions can also be symbolic over polynomial rings, and the \texttt{density} function allows one to calculate the density of various small patterns in a blowup. This example shows that in a large quasi-random graph containing a clique with relative size $x$ and edges between this clique and the rest with probability $p$, the density of induced cherries is $3px(x - 1)(3xp-2x-p)$.
\begin{lstlisting}[language=Python]
R.<x,p> = QQ[]
cherry = G(3, edges=[[0,1],[1,2]])
test_constr = G.blowup_construction(4, 
                                    [x,1-x], 
                                    edges={(0,0): 1, (0, 1):p}
                                   )
print(factor(test_constr.density(cherry)))
\end{lstlisting}

\paragraph{Certificates and problem files} The optimizer can write a certificate to a file and later verify it. Verification requires the same theory state (including exclusions) as used during generation. Setting the \texttt{file} parameter creates and writes the certificate to the specified file, if it exists, it overrides the file contents. This file can be later verified with the \texttt{Theory.verify} function.
\begin{lstlisting}[language=Python]
G.reset()
G.exclude(G(3))
G.optimize(G(2), 3, file="test")
G.verify("test")
\end{lstlisting}

To solve an SDP on a different machine or with a different solver, the SDP problem instance can be exported using \texttt{Theory.external\_optimize}. This expects the same optimization problem as \texttt{Theory.optimize}, but with the compulsory \texttt{file} parameter, to save the SDP problem instance. The exported SDP problem is stored in the SDPA sparse format.
\begin{lstlisting}[language=Python]
G.reset()
G.external_optimize(G(4, edges=[[0,1]]), 5, file="problem")
\end{lstlisting}

\section{Examples}\label{sec:examples}
This section contains scripts written for various results using the \texttt{FlagAlgebraToolbox}. The notebooks used to create these results and the corresponding certificates are available either in the ancillary folder of the arXiv version of the papers, or in the repository \url{https://github.com/bodnalev/supplementary_files}.

\subsection{Graph edge inducibility}

The following script, for a given $\kappa, \ell$ pair, tries to find the maximum density of $\kappa$ vertices inducing exactly $\ell$ edges a graph can have. This result appeared in \cite{bodnár2025turandensitytight5cycle}. To suppress the logs, the global \texttt{Theory.printlevel} disables all info and the optimization is called with a suitable construction and target size for each case. 

\begin{lstlisting}[language=Python]
G = GraphTheory
G.reset()
G.printlevel(0)
def strs(k, l):
    return sum([xx for xx in G.generate(k) if len(xx.edges)==l])

#lambda(3, 1)
constr = G.blowup_construction(
    5, 2, edges=[[0,0],[1,1]]
)
bound = G.optimize(
    strs(3, 1), 5, exact=True, 
    file="certificates/stats31", construction=constr
)
print("{}<=lambda(3, 1)<={}".format(
    constr.density(strs(3,1)), bound)
     )

#lambda(4, 1)
constr = G.blowup_construction(
    7, 5, edges=[[0,0],[1,1],[2,2],[3,3],[4,4]]
)
bound = G.optimize(
    strs(4, 1), 7, exact=True, denom=2**20, 
    file="certificates/stats41", construction=constr
)
print("{}<=lambda(4, 1)<={}".format(
    constr.density(strs(4,1)), bound)
     )

#lambda(4, 2)
constr = G.blowup_construction(
    7, 6, edges=[[0,1],[1,2],[2,0],[3,4],[4,5],[5,3]]
)
bound = G.optimize(
    strs(4, 2), 7, exact=True, denom=2**20, 
    file="certificates/stats42", construction=constr
)
print("{}<=lambda(4, 2)<={}".format(
    constr.density(strs(4,2)), bound)
     )


#lambda(4, 3)
constr = G.blowup_construction(
    6, 2, edges=[[0,1]]
)
bound = G.optimize(
    strs(4, 3), 6, exact=True, slack_threshold=1e-5, 
    denom=2**20, kernel_denom=2**20, 
    file = "certificates/stats43", construction=constr
)
print("{}<=lambda(4, 3)<={}".format(
    constr.density(strs(4,3)), bound)
     )
\end{lstlisting}
\begin{lstlisting}
3/4<=lambda(3, 1)<=3/4
72/125<=lambda(4, 1)<=72/125
1/2<=lambda(4, 2)<=1/2
1/2<=lambda(4, 3)<=1/2
\end{lstlisting}
Note that the paper considered all $\kappa \leq 7$ values and obtained other exact and non-exact results.

\subsection{The Turán density of the tight $5$-cycle minus one edge} This example script is from the paper \cite{BodnarPikhurko25ind}, where we determine the maximum edge density in $3$-graphs without $C_5^-$, which is the $3$-uniform graph forming a tight cycle on $5$ vertices with an edge missing. A super-saturation argument shows that the density of $K_4^-$ is also $0$, hence we can exclude both $K_4^-$ and $C_5^-$. The optimum is attained at a recursive blowup of $K_{3, 3, 3}$.

In the first script, we establish a weak upper bound with a basic usage of the software. 

\begin{lstlisting}[language=Python]
TGp = ThreeGraphTheory
k4m = TGp.pattern(4, edges=[[0,1,2],[0,1,3],[0,2,3]])
c5m = TGp.pattern(5, edges=[[0,1,2],[1,2,3],[2,3,4],[3,4,0]])
TGp.exclude([k4m, c5m])
edge = TGp(3, edges=[[0,1,2]])
a31 = TGp.optimize(edge, 7,  exact=True, 
                   file="certificates/proposition_3_1", 
                   denom=1024*3125, printlevel=0
                  )
print("Result is {} ~= {}".format(a31, a31.n()))
\end{lstlisting}
\begin{lstlisting}
Result is 401181/1600000 ~= 0.250738125000000
\end{lstlisting}

The second cell still works in the \texttt{TGp} theory, and first defines the quantity $F^E_{2, 2, 2}$ from the paper, then establishes that any $3$-graph without $K_4^-$ and $C_5^-$, and edge density close to optimum must have many $\llbracket F^E_{2, 2, 2} \rrbracket$ copies. This rough bound can be translated to a partitioning of the vertices into three parts, where the density of the $3$-edges between the parts is at least $0.194$. The idea is that these three parts closely approximate the top level structure of the optimal construction.

\begin{lstlisting}[language=Python]
b32 = 1/4 - 1/100000
f222 = TGp.pattern(6, ftype=[0,1,2], 
                   edges=[[0,1,2],[3,4,5],[0,1,5],[0,2,4],[1,2,3]])
k222 = f222.project()
gamma = TGp.optimize(k222, 7, maximize=False, positives=[edge - b32], 
                     exact=True, file="certificates/proposition_3_2", 
                     denom=1024*1024, printlevel=0
                    )
a32 = gamma / a31
print("Result is {} ~= {}".format(a32, a32.n()))
\end{lstlisting}
\begin{lstlisting}
Result is 1701468433/8763932672 ~= 0.194144398032181
\end{lstlisting}

The final script works in the theory where this three partition is specified by an additional binary relation. We combine the previous theory \texttt{TGp}, with a theory \texttt{P}, where there is only one binary relation, which is guaranteed to partition the vertices into three parts, by a suitable exclusion of small patterns.

The final calculation is performed in this combined theory. We wish to prove that any almost optimal $3$-graph, satisfying various sensible conditions and the condition from the previous calculation, must have the $3$-edge structure between parts matching the optimal construction. Essentially the extra binary relation allows us to mask out the recursive detail inside the parts, and only run a calculation arguing about what happens between. The rounding is aided by specifying the construction. For more details see the article \cite{bodnár2025turandensitytight5cycle}.

\begin{lstlisting}[language=Python]
# setting up partitioned theory
P = Theory("3Partition", relation_name="part", arity=2, is_ordered=False)
P.exclude([
    P(3, part=[[0,1]]), 
    P(4, part=[[0,1],[0,2],[0,3],[1,2],[1,3],[2,3]])
])
CTGp = combine("3PartitionNoC5m", TGp, P)

# edge with (C)orrect partition
C = CTGp(3, edges=[[0,1,2]],part=[[0,1],[0,2],[1,2]])
# edge with (C)orrect partition (p)ointed
Cp = CTGp(3, edges=[[0,1,2]], part=[[0,1],[0,2],[1,2]], ftype=[0])

# edge with (B)ad partition
B = CTGp(3, edges=[[0,1,2]], part=[[0,1],[1,2]])

# edge with (B)ad partition (p)ointed 
Bp = CTGp(3, edges=[[0,1,2]], part=[[0,1],[1,2]], ftype=[0])

# (M)issing edge with good partition
M = CTGp(3, edges=[], part=[[0,1],[0,2],[1,2]])

# at each point, good edges are more than bad edges 
b33 = 19/100
assums = [Cp - Bp/2, C - b33]

# optimal construction and its derivatives
symbolic_constr = \
CTGp.blowup_construction(6, 
                         ["X0", "X1", "X2"], 
                         edges=[[0, 1, 2]], 
                         part=[[0, 1], [0, 2], [1, 2]]
                        ).set_sum()
ders = symbolic_constr.derivatives([1/3, 1/3])

# bad is less than missing, proven by (B)ad minus (M)issing is at most 0.
CTGp.optimize(B - M*(99/100), 6, positives=assums, 
              exact=True, construction=ders, 
              file="certificates/proposition_3_3", 
              denom=1024*16, slack_threshold=1e-6, 
              kernel_denom=2**20, printlevel=0
             )
\end{lstlisting}

\section{Future work}\label{sec:future}

The package is actively used in research. Much of the functionality is custom-made for a specific problem and due to the lack of generalization, is not included in the full package. These provide possible functions we plan to integrate.  

\begin{itemize}
    \item \textbf{Plotting} A custom code plotting graph flags and graph flag algebra elements was valuable in multiple projects. Due to only supporting one theory, it is not included in the latest version of the package. The structures therefore only have a text-based description.
    \item \textbf{Stability} Perfect stability in optimization problems can often be automatically detected. The method was originally described in \cite{PikhurkoSliacanTyros19}. Ad-hoc functions were developed for this task in recent papers \cite{BodnarPikhurko25ind, bodnár2025semiinducibility4vertexgraphs}, but it was only tested for graphs.
    \item \textbf{Complex constructions} The current blow-up construction function does not support recursive constructions. In addition, the quasi-random relations are always independent. A potential extension would allow any joint probability distribution between different relations.
    \item \textbf{Interpretations} Mapping between different theories would be a valuable addition to the package, as described in \cite[Section 2.3]{Razborov07}. Currently this is only supported by custom functions, but a solid general functionality would be good for this.
    \item \textbf{Semidefinite assumptions} Currently, during optimization, the positivity assumptions with types are translated to linear assumptions, by multiplying them with all possible typed flags and then averaging the result. However, this could be extended to multiplication by sum-of-squares expressions, allowing a stronger usage of these constraints.
    \item \textbf{Fractional hierarchy} Currently, the SDP optimization collects and reasons over all flags up to a specified size. Unfortunately, the number of flags as the size increases is often super-polynomial. Allowing the reasoning over a small suitable fraction of flags with larger sizes would potentially allow better results, with less computation. 
\end{itemize}

\section{Availability, versioning}\label{sec:availability}

\texttt{FlagAlgebraToolbox} is distributed as a fork of the SageMath project at \url{https://github.com/bodnalev/sage}. Prebuilt images are distributed as \texttt{setempty/flag-algebra-toolbox} Docker images, and the launcher described in Subsection~\ref{subsec:installation} is available at \url{https://gist.github.com/bodnalev/900b7f8c1b4dcf0e69e232be361c04ff}. The code is released under the same open-source, GPL version 3 or later compatible license as SageMath. In addition to the standard components of SageMath, the software includes CSDP, a semidefinite programming solver \cite{Borchers01011999}, and Bliss, for calculating graph canonical labelings \cite{JunttilaKaski:ALENEX2007}.

This note documents \texttt{FlagAlgebraToolbox} version \texttt{0.3}. The corresponding release-specific Docker image is \texttt{setempty/flag-algebra-toolbox:0.3}.

If you use \texttt{FlagAlgebraToolbox} in academic work, please cite this paper and include the release tag used, as in the example below.
\begin{lstlisting}
@misc{flagalgebratoolbox,
  author = {Levente Bodn{\'a}r},
  title = {{FlagAlgebraToolbox}: Flag algebra computations in SageMath},
  year = {2026},
  note = {FlagAlgebraToolbox version 0.3},
  eprint = {2601.06590},
  archivePrefix = {arXiv},
  primaryClass = {math.CO},
  url = {https://arxiv.org/abs/2601.06590},
}
\end{lstlisting}

\section*{Acknowledgements}

This project was supported by the ERC Advanced Grant 101020255.

%\bibliography{bibexport}
\printbibliography

@article{BodnarPikhurko25ind,
  author = {Bodn{\'a}r, Levente and Pikhurko, Oleg},
  title = {Some exact inducibility-type results for graphs via flag algebras},
  journal = {The Electronic Journal of Combinatorics},
  volume = {33},
  number = {2},
  pages = {P2.39},
  year = {2026},
  doi = {10.37236/14612},
}

@article{PikhurkoSliacanTyros19,
  author = {Pikhurko, Oleg and Slia{\v{c}}an, Jakub and Tyros, Konstantinos},
  title = {Strong forms of stability from flag algebra calculations},
  journal = {Journal of Combinatorial Theory, Series B},
  volume = {135},
  pages = {129--178},
  year = {2019},
  doi = {10.1016/j.jctb.2018.08.001},
}

@article{Razborov07,
  author = {Razborov, Alexander A.},
  title = {Flag algebras},
  journal = {The Journal of Symbolic Logic},
  volume = {72},
  number = {4},
  pages = {1239--1282},
  year = {2007},
  doi = {10.2178/jsl/1203350785},
}

@misc{bodnár2025turandensitytight5cycle,
  author = {Bodn{\'a}r, Levente and Le{\'o}n, Jared and Liu, Xizhi and Pikhurko, Oleg},
  title = {The {Tur\'an} density of the tight 5-cycle minus one edge},
  year = {2024},
  eprint = {2412.21011},
  archivePrefix = {arXiv},
  primaryClass = {math.CO},
}

@article{bodnár2025semiinducibility4vertexgraphs,
  author = {Bodn{\'a}r, Levente and Pikhurko, Oleg},
  title = {Semi-inducibility of 4-vertex graphs},
  journal = {Discrete Applied Mathematics},
  volume = {392},
  pages = {303--323},
  year = {2026},
  doi = {10.1016/j.dam.2026.06.003},
}

@article{Borchers01011999,
  author = {Borchers, Brian},
  title = {{CSDP}, a {C} library for semidefinite programming},
  journal = {Optimization Methods and Software},
  volume = {11},
  number = {1--4},
  pages = {613--623},
  year = {1999},
  doi = {10.1080/10556789908805765},
}

@inproceedings{JunttilaKaski:ALENEX2007,
  author = {Junttila, Tommi and Kaski, Petteri},
  title = {Engineering an efficient canonical labeling tool for large and sparse graphs},
  booktitle = {Proceedings of the Ninth Workshop on Algorithm Engineering and Experiments (ALENEX 2007)},
  pages = {135--149},
  publisher = {SIAM},
  year = {2007},
  doi = {10.1137/1.9781611972870.13},
}

@article{flagmatic,
  author = {Falgas-Ravry, Victor and Vaughan, Emil R.},
  title = {{Tur\'an} {$H$}-densities for 3-graphs},
  journal = {The Electronic Journal of Combinatorics},
  volume = {19},
  number = {3},
  pages = {P40},
  year = {2012},
  doi = {10.37236/2733},
}

@manual{The_SageMath_Developers_SageMath_2026,
  author = {{The SageMath Developers}},
  title = {{SageMath}},
  version = {10.9},
  year = {2026},
  doi = {10.5281/zenodo.8042260},
  url = {https://www.sagemath.org},
}

@article{Mantel1907,
  author = {Mantel, Willem},
  title = {Vraagstuk {XXVIII}},
  journal = {Wiskundige Opgaven met de Oplossingen},
  volume = {10},
  number = {2},
  pages = {60--61},
  year = {1907},
}

@article{Even_Zohar_2014,
  author = {Even-Zohar, Chaim and Linial, Nati},
  title = {A note on the inducibility of 4-vertex graphs},
  journal = {Graphs and Combinatorics},
  volume = {31},
  number = {5},
  pages = {1367--1380},
  year = {2015},
  doi = {10.1007/s00373-014-1475-4},
}

\end{document}